\documentclass[noamsfonts]{amsproc}
\usepackage{amssymb}
\title[The Leading Coefficients of Certain Kazhdan-Lusztig  Polynomials]{The Leading
Coefficient of Certain Kazhdan-Lusztig \\  Polynomials of the
Permutation Group $\mathfrak S_n$}

\author[N. Xi]{Nanhua Xi}

\address{N.X.\\
Institute of Mathematics\\
Chinese Academy of Sciences \\
Beijing 100080\\
CHINA }
\email{nanhua@mail2.math.ac.cn}

\begin{document}

\begin{abstract}
In this paper we  show that the leading coefficient $\mu(y,w)$ of
certain Kazhdan-Lusztig polynomials $P_{y,w}$ of the permutation
group $\mathfrak S_n$ of 1,2,...,n are  not greater than 1. More
precisely, we show that the leading coefficients $\mu(y,w)$ are not
greater than 1 whenever $a(y)< a(w)$, where $a: \mathfrak
S_n\to\mathbf N$ is the function defined in [L1]. See 1.5 for a
simple interpretation of the function $a$.
\end{abstract}

\maketitle

\def\Cal{\mathcal}
\def\bold{\mathbf}
\def\ca{\mathcal A}
\def\cdz{\mathcal D_0}
\def\cd{\mathcal D}
\def\cdo{\mathcal D_1}
\def\bold{\mathbf}
\def\l{\lambda}
\def\le{\leq}

The Kazhdan-Lusztig polynomials of a Coxeter group $(W,S)$ play a
central role in Kazhdan-Lusztig theory. From the seminal work [KL]
one sees that the leading coefficient $\mu(y,w)$ of some
Kazhdan-Lusztig polynomials $P_{y,w}$ are of great importance. It
is certainly interesting to understand the leading coefficients.
We have two natural questions, (a) determine the leading
coefficients $\mu(y,w)$, (b) for a given element $w$ in the
Coxeter group, determine all $y$ such that $\mu(y,w)$ or
$\mu(w,y)$ is non-zero. By now, for both questions we do not have
much understanding except for some special cases, see for example
[L4].

\def\gsn{\mathfrak S_n}
\def\gs{\mathbf S}

In this paper we are concerned with question (a) for the
permutation group  $\gsn$ of 1,2,...,n. We shall show that for
$y,w\in\gsn$, if $y\le w$ (Bruhat order) and $a(y)< a(w)$, then
the leading coefficient $\mu(y,w)$ is less than or equal to 1, see
Theorem 1.4, here the function $a: \gsn\to \mathbf N$ is defined
by Lusztig in [L1], see 1.5 for a simple interpretation of the
function for $\gsn$. When $a(y)\ge a(w)$,  McLarnan and Warrington
have showed that $\mu(y,w)$ can be greater than 1, see [MW]. The
work [MW] gives a negative answer to the (0,1)-conjecture in [GM,
p.52].

In this paper we also give some discussion for the set of all
$\mu(y,w),\ y,w\in\gsn$, see section  1.2 and 3.2.

\section{Preliminaries}

In this section we collect some terms and basic facts, we refer to
[KL] for more details.

\medskip

 \noindent{\bf 1.1.}  Let $(W,S)$ ($S$ the
set of simple reflections) be a Coxeter group.  We shall denote the
length function of $W$ by $l$ and use $\leq$ for the Bruhat order
on $W$.

Let $H$ be the  Hecke algebra  of $(W,S)$ over $\Cal A=\bold
Z[q^{\frac 12},q^{-\frac 12}]$ $(q$ an indeterminate) with
parameter $q$. Let  $\{T_w\}_{w\in W}$ be its standard basis.
There is a unique ring homomorphism \ $\bar{\cdot}: H\to H$ such
that $\overline{T_w}=T^{-1}_{w^{-1}}$ and $\bar q^{\frac
12}=q^{-\frac 12}$.

In [KL1], Kazhdan and Lusztig proved that for each $w$ in $W$,
there exists a unique element $C_w$ in the Hecke algebra $H$ such
that (1) $C_w=\bar C_w,$ (2) $C_{w}=q^{-\frac{l(w)}2}\sum_{y\le
w}P_{y,w}T_y,$ where $P_{y,w}$ is a polynomial in $q$ of degree
less than or equal to $\frac 12(l(w)-l(y)-1)$ for $y< w$ and
$P_{w,w}=1$.

\medskip

The basis $\{C_w\}_{w\in W}$ of $H$ is the Kazhdan-Lusztig basis.
The polynomials $P_{y,w}$ are the Kazhdan-Lusztig polynomials. If
$y< w$, we have $$P_{y,w}=\mu(y,w)q^{\frac
12(l(w)-l(y)-1)}+\text{lower degree terms.}$$  We write $y\prec w$
if $\mu(y,w)\ne 0$. 

\medskip

\def\ll{\underset {L}{\leq}}
\def\rl{\underset {R}{\leq}}
\def\lrl{\underset {LR}{\leq}}
\def\llr{\lrl}
\def\el{\underset {L}{\sim}}
\def\er{\underset {R}{\sim}}
\def\elr{\underset {LR}{\sim}}
\def\ds{\displaystyle\sum}

 \noindent{\bf 1.2.} Assume that $(W,S)$ is an affine Weyl group or a Weyl
group. The
 following formula is due to Springer (see [Sp]),

\medskip

 \noindent(a)\qquad $\mu(y,x)=\displaystyle\sum_{d\in\cdz}\delta_{y^{-1},x,d}+
 \displaystyle\sum_{f\in\cdo}\gamma_{y^{-1},x,f}\pi(f)=\displaystyle\sum_{d\in\cdz}\delta_{y,x^{-1},d}+
 \displaystyle\sum_{f\in\cdo}\gamma_{y,x^{-1},f}\pi(f).$

 \medskip

 We need explain the notations.
Write $$C_xC_y=\sum_{z\in W}h_{x,y,z}C_z,\qquad h_{x,y,z}\in
\mathcal A
 =\bold Z[q^{\frac 12},q^{-\frac 12}].$$
 Following Lusztig ([L1]), we define
 $$a(z)={\rm min}\{i\in\bold N\ |\ q^{-\frac i2}h_{x,y,z}\in\bold Z[q^{-\frac
 12}]{\rm\ for \ all\ }x,y\in W\}.$$ If for any $i$,
 $q^{-\frac i2}h_{x,y,z}\not\in\bold Z[q^{-\frac
 12}]{\rm\ for \ some\ }x,y\in W$, we set $a(z)=\infty.$
 It is not clear that whether there exists a Coxeter group
 $(W,S)$ such that $a(z)=\infty$ for some $z\in W$.

 From now on, we assume that the function $a:\ W\to\bold N$ is
 bounded and $(W,S)$ is crystallographic. In this case, the function $a$ is
 constant on two-sided cells of $W$ (see [L1]). Obviously, when $W$ is
 finite, the function $a$ is bounded. Lusztig showed that the
 function $a$ is  bounded for all affine Weyl groups (see [L1]).
 Following Lusztig and Springer,
 we define $\delta_{x,y,z}$ and $\gamma_{x,y,z}$ by the following formula,
 $$h_{x,y,z}=\gamma_{x,y,z}q^{\frac {a(z)}2}+\delta_{x,y,z}q^{\frac
 {a(z)-1}2}+
 {\rm\ lower\ degree\ terms}.$$

 Springer showed that $l(z)\ge a(z)$ (see [L2]). Let $\delta(z)$ be the
 degree of $P_{e,z}$, where $e$ is the neutral element of $W$.
 Then actually one has $l(z)-a(z)-2\delta(z)\ge 0$ (see [L2]). Set
 $$\cd_i=\{z\in W\ |\ l(z)-a(z)-2\delta(z)=i\}.$$The number $\pi(z)$ is
 defined by $P_{e,z}=\pi(z)q^{\delta(z)}+
 {\rm\ lower\ degree\ terms}.$

 The elements of $\cdz$ are involutions, called distinguished involutions of
 $(W,S)$ (see [L2]). Moreover, in a Weyl group or an affine Weyl
 group, each left cell (resp. right cell) contains exactly one element  of $\cdz$, see [L2].
 When $W$ is a Weyl group, Lusztig has showed that an element $z$ is in $\cdo$
if and only if there exists some $d$ in $\cdz$ such that $z\elr d$
and $\mu(z,d)\ne 0$ or $\mu(d,z)\ne 0$, see [L3].
 For $\gsn$, all involutions are distinguished.

\medskip

 Example: Let $S'$ be a subset of $S$ such that the subgroup $W'$
 of $W$ generated by $S'$ is finite. Then the longest element of
 $W'$ is in $\cdz$.

 \medskip

\medskip

\noindent{\bf 1.3.} We refer to [KL] for the definition of the
preorders $\ll,\ \rl,\lrl$ and of the equivalence relations $\el,\
\er,\ \elr$ on $W$. The corresponding equivalence classes are
called { left cells, right cells, two-sided cells} of $W$,
respectively. The preorder $\ll$ (resp. $\rl;\lrl$) induces a
partial order on the set of left (resp. right; two-sided) cells of
$W$, denoted again by $\ll$ (resp. $\rl;\lrl$). For a Weyl group
or an affine Weyl
 group, Springer showed
the following results (a) and (b) (see [Sp])

\medskip

\noindent (a) Assume that $\mu(y,w)$ or $\mu(w,y)$ is nonzero,
then
 $y\ll w$ and $y\rl w$ if $a(y)<a(w)$, and $y\el w$ or $y\er w$ if
 $a(y)=a(w)$.

\medskip

\noindent (b) If $\delta_{x,y,z}\ne 0$, then $z\el y$ or $z\er x$.
 (Note that  $h_{x,y,z}\ne 0$ implies that $a(z)\ge a(x)$ and
 $a(z)\ge a(y)$, see [L1].)

\medskip

For $w\in W$, set $L(w)=\{s\in S\ |\ sw\leq w\}$, \ $R(w)=\{s\in
S\ |\ ws\leq w\}.$ Then we have (see [KL])

\medskip

\noindent(c) $R(w)\subseteq R(y)$ if $y\ll w$. In particular,
$R(w)=R(y)$ if $y\el w$;

\medskip

\noindent(d) $L(w)\subseteq L(y)$ if $y\rl w$. In particular,
$L(w)= L(y)$ if $y\er w$.

\medskip

We shall need the star operation introduced in [KL].
 Let $s$ and $t$ be in $S$
such that $st$
 has order 3, i.e. $sts=tst$. Define
$$D_L(s,t)=\{w\in W\  |\   L(w)\cap\{s,t\}\text{ has exactly one
element}\},$$
$$D_R(s,t)=\{w\in W\ |\ R(w)\cap\{s,t\}\text{ has
exactly one element}\}.$$

\medskip

If $w$ is in $D_L(s,t)$, then $\{sw,tw\}$ contains exactly one
element in $D_L(s,t)$, denoted by ${}^*w$, here $*=\{s,t\}$. The
map: $\varphi: \ D_L(s,t)\to D_L(s,t)$, $w\to {}^*w$, is an
involution and is called a {  left star operation}. Similarly if
$w\in D_R(s,t)$ we can define the {  right star operation} $\psi:\ D_R(s,t)\to D_R(s,t),\ 
w\to w^*=\{ws,wt\}\cap D_R(s,t)$, where $*=\{s,t\}$.
The following result is proved in [KL].

\medskip

\noindent(e) Let  $s$ and $t$ be in $S$ such that $st$
 has order 3. Assume that $y\prec w$. If both $y^*$ and $w^*$ (resp.
 ${}^*y$ and ${}^*w$) are well defined , then either (1) $y^*\prec w^*$ and $\mu(y,w)=\mu(y^*,w^*)$, or
 (2) $w^*\prec y^*$ and $\mu(y,w)=\mu(w^*,y^*)=1$ (resp. either (1)
 ${}^*y\prec{}^*w$ and $\mu(y,w)=\mu({}^*y,{}^*w)$, or (2)${}^*w\prec{}^*y$ and
 $\mu(y,w)=\mu({}^*w,{}^*y)=1.$)

\medskip

The main result of this paper is the following.

 \noindent{\bf 1.4. Theorem.}  Let $y,w\in\gsn$. If $y< w$  and
 $a(y)<
a(w)$, then the leading coefficient $\mu(y,w)$ is less than or
equal to 1, here the function $a: \gsn\to \mathbf N$ is defined in
1.2.

\medskip

{\bf Remark.} This result is also true for an affine Weyl group of
type $\tilde A_{n-1}$.

\medskip

  \noindent{\bf 1.5.}
 For $W=\gsn$, there is a nice description for cells of
$W$, called Robinson-Schensted rule. We recall the rule since it
also provide a simple way to compute the function
$a:\gsn\to\mathbf N$ defined in subsection 1.2. As usual, the set
$S$ consists of transpositions $s_i=(i,i+1)$, $i=1,2,...,n-1.$

A partition $\l$ of $n$ is visualized as a Young diagram $F_\l$
with $n$ boxes such that the lengths of the columns of the diagram
are given by the components of $\l$. A Young tableau of shape $\l$
is a labelling of the boxes of $F_\l$ with integers 1,2,...,n. The
tableau is called standard if the numbers labelling its boxes are
increasing in the columns from the top to the bottom and are
increasing in the row from left to right. The Robinson-Schensted
rule gives a bijection $\theta: w\to (P(w),Q(w)$ from the
permutation group $\gsn$ to the set of pairs of standard Young
tableau on $\{1,2,...,n\}$. Let $w(i)=j_i$ for $i=1,2,...,n$. The
standard Young tableau $P(w)$ is defined as follows.

Let $T_0$ be the empty Young tableau. Assume that we have defined
a Young tableau $T_{i-1}$ which has totally $i-1$ boxes and the
boxes are labelled with integers $j_1,j_2,...,j_{i-1}$. Now we
define the Young tableau $T_i$. If $j_i$ is greater than every
numbers in the first row of $T_{i-1}$, then $T_i$ is obtained from
$T_{i-1}$ by adding $j_i$ to the end of  the first row of
$T_{i-1}$. If not, then let $z$ be the smallest number of the
first row of $T_{i-1}$ greater than $j_i$. Let $T'_{i-1}$ be the
Young tableau by removing the first row of $T_{i-1}$. Then $T_i$
is the Young tableau obtained from $T_{i-1}$ by replacing $z$ by
$j_i$ and by inserting $z$ to the tableau $T'_{i-1}$ as the same
way. It is known that $P(w)$ is indeed a standard Young tableau.
We set $Q(w)=P(w^{-1})$. Let $y,w$ be in $\gsn$, then we have (see
[BV])

\medskip

\noindent(a) $y\el w$ if and only if $Q(y)=Q(w)$,

\medskip

\noindent(b) $y\er w$ if and only if $P(y)=P(w)$,

\medskip

\noindent(c) $y\elr w$ if and only if $P(y$ and $P(w)$ have the
same shape.

\medskip
\def\l{\lambda}

Assume that $P(w)$ is a Young tableau of shape
$\l=(\l_1\ge\l_2\ge\cdots\ge\l_k)$ (a partition of $n$). we define
 $a(w)$ to be $\frac 12 \sum_{i=1}^k\l_i(\l_i-1)$. Then
$a:\gsn\to\mathbf N$ is exactly the $a$-function defined in 1.2.

Note that $y\llr w$ if and only if $\l(y)\le \l(w)$, where $\l(y)$
and $\l(w)$ are the partitions corresponding to $P(y)$ and $P(w)$
respectively. (We say $\lambda=(a_1,a_2,...)\le\mu=(b_1,b_2,...)$
if $a_1+\cdots+a_k\le b_1+\cdots+b_k$ for all $k=1,2,...$.)

\medskip

\section{A description of the set of all $\mu(y,w)$ for $\gsn$}

In this section we give a description to the set of all
$\mu(y,w),\ y,w\in\gsn$, in terms of $\delta_{x,y,z},\
x,y,z\in\gsn$. More precisely, we have

\medskip

\noindent{\bf Proposition 2.1.} For $\gsn$, the set
$\{\delta_{x,y,z}|x,y,z\in\gsn\}$ is equal to the set\break 
$\{\mu(y,w)|y,w\in\gsn\}$.

\medskip

{\it Proof:} Clearly both sets contain 0. Assume that $\mu(y,w)$
is non zero. By Springer's formula 1.2(a), we have
$\delta_{y^{-1},w,d}\ne 0$ for some $d\in\cdz$ or
$\gamma_{y^{-1},w,u}\ne 0$ for some $u\in \cdo$. It is no harm to
assume that $a(w)\ge a(y)$. If $a(w)>a(y)$, then
$\gamma_{y^{-1},w,u}= 0$ for all $u$ in $\gsn$ (see [L1, 6.3(a)]),
and $\delta_{y^{-1},w,d}\ne 0$ implies that $d$ and $w$ are in the
same left cell, cf 1.2 (c) and [L1, 6.3(b)]. Therefore
$\mu(y,w)=\delta_{y^{-1},w,d(w)}$ in this case. (We use $d(x)$ for
the unique element in $\cdz$ which is in the left cell containing
$x$.) If $a(w)=a(y)$, by 1.3 (a), then $w$ and $y$ are in the same
left cell or in the same right cell. It is no harm to assume that
$y,w$ are in the same left cell since
$\mu(y,w)=\mu(y^{-1},w^{-1})$. Then $y,w$ are in different right
cell (cf [KL,\S5]. This implies that $\gamma_{y^{-1},w,u}= 0$ for
all $u\in\gsn$ (see [L1, 6.3(a)]. If $\delta_{y^{-1},w,d}\ne 0$,
then $d=d(w)$. Thus in this case we still have
$\mu(y,w)=\delta_{y^{-1},w,d(w)}$.

Conversely, assume that $\delta_{x,y,z}\ne 0$. It is no harm to
assume that $a(y)\ge a(x)$. Then $y$ and $z$ are in the same left
cell, cf 1.2 (c) and [L1, 6.3(b)]. Thus we can find a sequence of
right star operations $\varphi_k,...,\varphi_{1}$ such that
$\varphi_k\cdots\varphi_{1}(z)=d$ is in $\cdz$, see [KL, \S5]. By
a similar argument as for [X, Prop. 1.4.4(b)], we see that
$h_{x,y,z}=h_{x,y',d},$ here $y'=\varphi_k\cdots\varphi_{1}(y)$.
In particular we have $\delta_{x,y,z}=\delta_{x,y',d}.$
 If $a(y)>a(x)$,
by Springer's formula 1.2 (a), the value $\delta_{x,y',d}$ is equal
to $\mu(x^{-1},y')$ or $\mu(y',x^{-1})$. If $a(y)=a(x)$, then $d$
and $x$ are in the same right cell, $d$ and $y'$ are in the same
left cell. So $x,{y'}^{-1}$ are in the same right cell and are not
in the same left cell. Thus $\gamma_{x,y',u}=0$ for any $u$ (see
[L1, 6.3(a)]. In this case we also have
$\delta_{x,y',d}=\mu(x^{-1},y')$ or
$\delta_{x,y',d}=\mu(y',x^{-1})$.

The proposition is proved.

\medskip

In the proof we see the following

\noindent{\bf Corollary 2.2.}
 If $a(y)\ge a(x),$ then $\delta_{x,y,z}=\delta_{x,y',d}=\mu(x^{-1},y')$ or $\mu(y',x^{-1})$ for some $y'$ in
the right cell containing $y$ and some $d$ in $\cdz$.

\medskip

When $y$ and $w$ are in the same two-sided cell we have the
following result.

\noindent{\bf Proposition 2.3.} The set $\{\pi(f)|f\in\cdo\}$ is
the same as the set $\{\mu(x,y)|x\elr y{\text{ and }} \mu(x,y)\ne
0\}$.

{\it Proof:} If $\mu(x,y)\ne 0$ and $x$ and $y$ are in the same
two-sided cell, by 1.3 (b) and [L1, 6.3],  $x$ and $y$ are in the
same left cell or in the same right cell. It is no harm to assume
that $x$ and $y$ are in the same left cell. By Springer's formula
1.2 (a), $\delta=\delta_{x,y^{-1},d(y^{-1})}\ne 0$ or
$\gamma=\gamma_{x,y^{-1},f}\ne 0$ for some $f\in\cdo$. Since $y$
and $w$ are not in the same right cell, we have $\delta=0$ and
$\gamma\ne 0$. Note that $\gamma_{x,y^{-1},f}\ne 0$ implies that
$f\er x$ and $f\el y^{-1}$, see [L1, 6.3]. Such $f$ is unique.
Moreover, using [X, Theorem 1.4.5], [KL, \S5] and [L2, 1.4(a)],
$\gamma_{x,y^{-1},f}=1$ in this case. Therefore
$\mu(x,y)=\gamma_{x,y^{-1},f}\pi(f)=\pi(f).$ The proposition is
proved.

\medskip

For simplifying statements, {\sl in the rest of this section we
set $\mu(y,w)=\mu(w,y)$ whenever one of $\mu(y,w)$, $\mu(w,y)$ is
defined.}  Note that we have
$\mu(f,d(f))=\gamma_{f,d(f),f}\pi(f)=\pi(f)$. Thus we have

\medskip

\noindent{\bf Corollary 2.4.} The set $\{\mu(x,y)|x\elr y{\text{
and }} \mu(x,y)\ne 0\}$ is the same as the set
$\{\mu(f,d)|f\in\cdo,\ d\in\cdz,{\text{ and }}f\elr d,\
\mu(f,d)\ne 0\}.$

\medskip

\noindent{\bf 2.5.} Let $f,x\in\gsn$. If $x\elr f$ and $h_{f,f,x}\ne
0$, then $x\el f$ and $x\er f$. Thus we must have $x=f$, since in a
two-sided cell of $\gsn$, each left cell and each right cell have
exactly one common element. Let $M_f$ be the irreducible module of
$H\otimes_A\mathbf Q(q^{\frac 12})$ provided by the left cell of $\gsn$
containing $f$ (cf, [KL]). Then we have
$$h_{f,f,f}=tr(C_f,M_f).$$ By [X, Prop. 1.4.5 (b)], we see that
$h_{f,f,f}=h_{f,d,d}$, where $d\in\cdz$ and $d\er f$. Thus
$\delta_{f,f,f}=\mu(f,d)$. Therefore we also have
$\{\mu(x,y)|x\elr y\}=\{\delta_{f,f,f}|f\in\gsn\}$ (we set
$\mu(x,x)=0$).

\medskip

\section{Proof of  Theorem 1.4}

\medskip

\noindent{\bf 3.1.} Now we can prove the main result. Let $y\le w$
be elements in $\gsn$ such that $a(y)<a(w)$ and $\mu(y,w)\ne 0$.

By 1.3 (a), we have $y\ll w$ and $y\rl w$. Thus $L(y)$ is a subset
of $L(w)$. If there is a simple reflection $s$ in $L(w)-L(y)$,
then $P_{y,w}=P_{sy,w}$ and the condition $\mu(y,w)\ne 0$ forces
that $sy=w$ and $\mu(y,w)=1$. Now assume that $L(w)=L(y)$.
According to [KL, \S5], there exists some subset $I$ of
$\{1,2,...,n-1\}$ and a sequence $\varphi_i,...,\varphi_1$ of left
star operations such that $\varphi_i\cdots\varphi_1(w)\er w_I$.
(We use $W_I$ for the subgroup of $\gsn$ generated by all $s_i,\
i\in I$, and denote by $w_I$ the longest element of $W_I$.) For
$j=1,2,...,i,$ set $w_j=\varphi_j\cdots\varphi_1(w)$ and set
$w_0=w$.

Let $k$ be the maximal number among 1,2,...,$i$ such that
$y_j=\varphi_j\cdots\varphi_1(y)$ is well defined for all
$j=1,2,...,k$. Set $y_0=y$. If for some $1\le j\le k$, we have
$w_j\le y_j$ and $y_m\le w_m$ for $m=0,1,...,j-1$, using 1.3 (e)
we see that $1=\mu(w_j,y_j)=\mu(y_{j-1},w_{j-1})=\mu(y,w)$.

Now assume that $y_j\le w_j$ for all  $j=1,2,...k$. Using 1.3 (e)
we get $\mu(y_{k},w_{k})=\mu(y,w).$ We claim that $L(w_k)-L(y_k)$
is non-empty. Since
 $a(y_j)=a(y)<a(w)=a(w_j)$ for all $j=0,1,2,...,k,$  we
always have $L(y_j)\subseteq L(w_j)$ for all $j$. If $k< i$, the
set $L(w_k)-L(y_k)$ must be non-empty. Otherwise, we must have
$L(w_k)=L(y_k)$. Thus $y_{k+1}=\varphi_{k+1}(y_k)$ is well
defined. This contradicts the assumption on $k$. When $k=i$, we
have $L(w_k)=\{s_j | j\in I\}$. If $L(w_k)\subseteq L(y_k)$, then
we can find $y'$ such that $y_k=w_Iy'$ and $l(y_k)=l(w_I)+l(y')$.
Thus $a(y)=a(y_k)\ge a(w_I)=a(w_k)=a(w)$. A contradiction.
Therefore in this case we also have that $L(w_k)-L(y_k)$ is
non-empty.  We then must have $l(w_k)=l(y_k)+1$ and
$\mu(y_k,w_k)=\mu(y,w)=1$.

The theorem is proved.

\medskip

{\bf Remark.} Using the result [S, Lemma 18.3.2] and a similar argument as above one
can prove that Theorem 1.4 remains true for an affine Weyl group of
type $\tilde A_{n-1}$.

\medskip

\noindent{\bf Proposition 3.2.} For $\gsn,\ n\ge 4$, the set
$\{\mu(y,w)|y,w\in\gsn\}$ is equal to the set
$\{\mu(y,w)|y,w\in\gsn\text{ and }a(y)>a(w)\}$.

\medskip

{\it Proof:}  Clearly, when $n\ge 4$, 1 is in the set
$A=\{\mu(y,w)|y,w\in\gsn\text{ and }a(y)>a(w)\}$. By Theorem 1.4,
we only need to show that if $a(y)=a(w)$, then $\mu(y,w)$ is in
$A$. To do this we need some preparation.

\medskip

 Let $U=U_0$ be the subgroup of $\gsn$ generated by
$s_2,s_3,...,s_{n-1}$. For $i=1,2,...,n-1$, we set $U_i=s_i\cdots
s_2s_1U_0$ and $U'_i=U_0s_1s_2\cdots s_{i}$. Then $U_i\cap
U_j=\emptyset$ if $i\ne j$ and $\gsn=U_0\cup U_1\cup\cdots\cup
U_{n-1}$. Moreover if $w$ is in $U_0$, then the length of
$s_i\cdots s_2s_1 w$ is $l(w)+i$. The following assertions (a) and
(b) are easy to prove.

\medskip

\noindent(a) Let $y\in U_i$ and $w\in U_j$. If $y\le w$, then
$i\le j$.

\medskip

\noindent(b) Let $y=s_i\cdots s_2s_1 y_1$ and $w=s_i\cdots s_2s_1
w_1$ be elements in $U_i$. Then $P_{y,w}=P_{y_1,w_1}.$

For a subset $I$ of $\{1,2,...,n-1\}$, let $W_I$ be the subgroup
of $\gsn$ generated by $s_i,i\in I$ and $w_I$ the longest element
of $W_I$.

Let $I=I_1\cup I_2\cup\cdots\cup I_m$ be a subset of $\{1,2,...,
n-1\}$. Assume that $I_1=\{1,2,...,k_1\}$ and
$I_j=\{k_{j-1}+2,k_{j-1}+3,...,k_j\}$ for $j=2,3,...,m$. Then
$w_I=w_{I_1}w_{I_2}\cdots w_{I_m}$. Assume that $k_1+1\ge
k_2-k_1\ge k_3-k_2\ge \cdots\ge k_m-k_{m-1}\ge 1.$ The shape of
the Young tableau $P(w_I)$ is
$\l(w_I)=(k_1+1,k_2-k_1,...,k_m-k_{m-1},1,...,1)$.

According to 1.5 and  [KL, \S5], we have

\noindent(c) for each element $w$ in $\gsn$ we can find a unique
subset $I_\l$ of $\{1,2,...,n-1\}$  as above such that $w\elr w_I$
and can find a sequence $\varphi_i,...,\varphi_1$ of left star
operations such that $\varphi_i\cdots\varphi_1(w)\er w_I$.

Now we claim the following assertion.

\noindent(d) Let $I$ be as above and $w\er w_I$. Then $w$ is in $
U_{k_1}$.

We argue for (d). Let $j$ be such that $w$ is in $ U_j$. Clearly
we have $j\ge k_1$. If $j>k_1$, then $j=k_a$ for some $2\le a\le
m$ and $w=s_js_{j-1}\cdots s_2s_1w_1$ for some $w_1\in U$. Let
$u=s_{k_{a-1}+2}s_{k_{a-1}+3}\cdots s_jw$. Then $w\lrl u$ and
$s_{k_{a-1}+1}u\le u$. Clearly, for any $i$ in $ I_1\cup I_2\cup
\cdots\cup {I_{a-1}}\cup I_{a+1}\cup \cdots\cup I_m$, we have
$s_is_{k_{a-1}+2}s_{k_{a-1}+3}\cdots
s_j=s_{k_{a-1}+2}s_{k_{a-1}+3}\cdots s_js_i$. Thus, for these $i$
we have $s_iu\le u$. It is easy to check that for $i$ in
$I_a-\{k_{a-1}+2\}$, we have $s_iu\le u$. Let $$J=I_1\cup I_2\cup
\cdots I_{a-2}\cup{(I_{a-1}\cup\{k_{a-1}+1\})}\cup
(I_{a}-\{k_{a-1}+2\})\cup I_{a+1}\cup \cdots\cup I_m.$$ Then
$u=w_Ju_1$ for some $u_1$ with $l(u)=l(w_J)+l(u_1)$. Thus $a(u)\ge
a(w_J)=a(w_I)+1+(k_{a-1}-k_{a-2})-(k_a-k_{a-1})\ge a(w_I)+1$ (we
set $k_{-1}=-1$). This contradicts $a(u)\le a(w)=a(w_I)$.
Therefore (d) is true.

Now we can complete the proof of the proposition. Assume that
$a(y)=a(w)$. Then $y\el w$ or $y\er w$. It is no harm to assume
that $y\er w$ since $\mu(y,w)=\mu(y^{-1},w^{-1})$. Also we may
assume that $\mu(y,w)>1$. Since a left star operation sends a
right cell to a right cell with the same $a$-function value and it
keeps the value of the $\mu(y,w)$ (see 1.3(e)), it is now harm to
assume that $y\er w\er w_I$, where $I$ is as in (d). By (d), we
know that both $y$ and $w$ are in $U_{k_1}$. Set $y=s_{k_1}\cdots
s_2s_1y_1$ and $w=s_{k_1}\cdots s_2s_1w_1$. Then $y_1$ and $w_1$
are in $U$. By (b), we have $\mu(y,w)=\mu(y_1,w_1)$. If
$a(y_1)<a(w_1)$ or $a(y_1)>a(w_1)$, we are done. If
$a(y_1)=a(w_1)$,  continuing the above process or using induction
on $n$, we see that $\mu(y,w)$ is in $A$.

The proposition is proved.

\medskip

Here is a question. For $y\le w$ in $\gsn$ with $\mu(y,w)\ne 0$ and $a(y)>a(w)$, whether
we can find $y'$ and $w'$ in some $\mathfrak S_r$ such that $y'\elr w'$ and $\mu(y,w)=\mu(y',w')$.

\medskip

\noindent{\bf 3.3.} Using 1.5 (b) we get the following assertion.

\noindent(a) Let $w$ be in $\gsn$ and $I$ be as (d) of the proof
of Proposition 3.2. Then $w\er w_I$ if and only if the following
conditions are satisfied:

\medskip

\noindent (1) For $k_{i-1}+2\le f<g\le k_i$ (set $k_0=-1$), we
have $w^{-1}(f)>w^{-1}(g)$.

\medskip

\noindent(2) For each $1\le f\le n$, we have
$|w(I_i)\cap\{1,2,...,f\}|\ge |w(I_{i+1})\cap\{1,2,...,f\}|.$ (We
set $I_{m+j}=\{k_m+j\}$ for $1\le j\le n-k_m$.)

\medskip

\noindent(b) By means of (b) in the proof of Proposition 3.2 we see easily that
$\mu(y,w)\le 1$ if $y\le w$ and $a(w)=1$.

\bigskip

\noindent{\bf Acknowledgement:} Part of the work was done during
my visit to the Department of Mathematics, Hong Kong University of
Science and Technology. I am very grateful to the department for
hospitality and for financial support. I would like to thank Professor 
Y. Zhu for helpful discussions and for arranging the visit.


\end{document}